\documentclass{amsart} 
 \newtheorem{thm}{Theorem}[section]
 \newtheorem{cor}[thm]{Corollary}
 \newtheorem{lem}[thm]{Lemma}
 \newtheorem{prop}[thm]{Proposition}
 \theoremstyle{definition}
 
 \theoremstyle{remark}

 \numberwithin{equation}{section}

\newcommand{\proofend}{\hfill $\Box$ \\[.25in] }
\newcommand{\scal}[1]{\langle#1\rangle}

\newcommand{\C}{\mathbb C}
\newcommand{\R}{\mathbb R}

\newcommand{\rarrow}{\rightarrow}
\newcommand{\CA}{\mathcal{A}}
\newcommand{\M}{{\mathbb M}}
\newcommand{\va}{\varepsilon}
\newcommand{\zl}{\lambda}
\newcommand{\diag}{\operatorname{diag}}
\newcommand{\Z}{{\mathbb Z}}
\newcommand{\T}{{\mathbb T}}

\begin{document}
%
%
%
%
%
%
%
%
\title[Algebras of $AP$ functions: Hermite property]
{Algebras of Almost Periodic Functions with Bohr-Fourier Spectrum in a
Semigroup: Hermite Property and its Applications}

\author[Rodman]{Leiba Rodman}

\address{%
Department of Mathematics \\
College of William and Mary \\
Williamsburg, VA 23187-8795\\
USA} \email{lxrodm@math.wm.edu [LR], ilya@math.wm.edu [IS]}

\thanks{The research of both authors was partially supported by NSF grant DMS-0456625.}
\author[Spitkovsky]{Ilya M. Spitkovsky}
\subjclass{Primary 42A75; Secondary 46C20, 46J20, 47A68, 47B35}

\keywords{Almost periodic functions, Hermite rings, Wiener algebra,
matrix functions, factorization, Toeplitz corona}

\begin{abstract}
It is proved that the unital Banach algebra of almost periodic
functions of several variables with Bohr-Fourier spectrum in a given
additive semigroup is an Hermite ring. The same property holds for
the Wiener algebra of functions that in addition have absolutely
convergent Bohr-Fourier series. As applications of the Hermite
property of these algebras, we study factorizations of Wiener--Hopf
type of rectangular matrix functions and the Toeplitz corona problem
in the context of almost periodic functions of several variables.
\end{abstract}

\maketitle


\section{Introduction}
\setcounter{equation}{0}

We let $\R$ be the real field, $\R^k$ the vector space of real
$k$-dimensional vector columns, and $AP^k$ the algebra of complex
valued almost periodic functions of $k$ real variables $x=(x_1,
x_2,\ldots, x_k)$, i.e., the closed subalgebra of $L^\infty(\R^k)$
(with respect to the standard Lebesgue measure) generated by all the
functions $e_\lambda(x)= e^{i\langle \lambda, x \rangle}$. Here the
variable $x=(x_1, \cdots ,x_k)^T$ and the parameter
$\lambda=(\lambda_1, \cdots ,\lambda_k)^T$ are in $\R^k$, the
superscript $T$ denotes the transposition operation on a vector (or
a matrix), and
$$
\langle \lambda , x \rangle = \sum_{j=1}^k \lambda_jx_j
$$
is the standard inner product of $\lambda$ and $x$. The norm in
$AP^k$ will be denoted by $\|\cdot\|_{\infty}$.

The next proposition is standard (see, e.g.,  Section~1.1 in
\cite{Pan}).

\begin{prop}\label{pr:2.1}
$AP^k$ is a commutative unital $C^*$-algebra, and therefore can be
identified with the algebra $C({\mathcal B}_k)$ of complex valued
continuous functions on a certain compact Hausdorff topological
space ${\mathcal B}_k$. Moreover, $\R^k$ is dense in ${\mathcal
B}_k$.
\end{prop}

The space ${\mathcal B}_k$ is called the {\em Bohr compactification}
of $\R^k$.

For any $f\in AP^k$ its {\em Bohr-Fourier series} is defined by the
formal sum
\begin{equation}\label{mar281}
\sum_\lambda f_\lambda e^{i\langle \lambda, t\rangle},\quad
t\in\R^k,
\end{equation}
where
\[ 
f_\lambda = \lim_{T \rarrow \infty} \ \frac{1}{(2T)^k}
\int_{[-T,T]^k} e^{-i\langle\lambda, x\rangle } f(x)dx, \ \lambda
\in \R^k,\]
and the sum in (\ref{mar281}) is taken over the set $\sigma(f) =
\{\lambda \in \R^k\colon f_\lambda \neq 0\}$, called the {\em
Bohr-Fourier spectrum} of $f$. The Bohr-Fourier spectrum of every $f
\in AP^k$ is at most a countable set. The {\em Bohr mean} $M\{f\}$
of $f \in AP^k$ is given by $$M\{f\} := f_0 = \lim_{T \rarrow
\infty} \frac{1}{(2T)^k} \int_{[-T,T]^k} f(x)dx. $$

The {\em Wiener algebra} $APW^k$ is defined as the set of all $f \in
AP^k$ such that the Bohr-Fourier series (\ref{mar281}) of $f$
converges absolutely. The Wiener algebra is a Banach $*$-algebra
with respect to the {\em Wiener norm} $\|f\|_W=\sum_{\lambda\in
\R^k} |f_{\lambda}|$ (the multiplication in $APW^k$ is pointwise).
Note that $APW^k$ is dense in $AP^k$. For the general theory of
almost periodic functions of one and several variables we refer the
reader to the books \cite{C,Le,LeZh} and to Chapter~1 in \cite{Pan}.

Let $\Delta$ be a non-empty subset of $\R^k$. Denote
\[
AP^k_\Delta  =  \{f \in AP^k\colon \sigma(f) \subseteq \Delta\},
\quad APW^k_\Delta  =  \{f \in APW^k\colon \sigma(f) \subseteq
\Delta\}.
\]
If $\Delta$ is an additive subset of $\R^k$, then $AP^k_\Delta$
(resp. $APW^k_\Delta$) is a subalgebra of $AP^k$ (resp. $APW^k$),
which is unital if in addition $0\in \Delta$.

A subset $S$ of $\R^k$ is said to be a {\em halfspace} if it has the
following properties:
\begin{itemize}
\item[(i)] $\R^k=S\cup(-S)$;
\item[(ii)] $S\cap(-S)=\{0\}$;
\item[(iii)] if $\lambda,\mu \in S$ then $\lambda+\mu\in S$;
\item[(iv)] if $\lambda \in S$ and $\alpha$ is a nonnegative real number, then
$\alpha \lambda \in S$.
\end{itemize}
A standard example (used extensively in \cite{RSW01ot,RSW02jot}) of
a halfspace is given by the set \[ E_k=\{(\lambda_1, \cdots ,\lambda_k)^T \in
\R^k \colon \lambda_1=\lambda_2= \cdots =\lambda_{j-1}=0, \
\lambda_j \neq 0 \ \Rightarrow \lambda_j>0 \}.
\]

Throughout the paper every ring is assumed to be commutative and
with unity, which will be denoted $e$. Let ${\mathcal R}$ be a ring.
An ordered $n$-tuple $(a_1, \ldots ,a_n)$ of elements of ${\mathcal
R}$ is said to be {\em unimodular}, if there exist $b_1, \cdots ,b_n
\in {\mathcal R}$ such that $a_1b_1+ \cdots +a_nb_n=e$. A ring
${\mathcal R}$ is called an {\em Hermite ring} if every unimodular
row can be complemented; in other words, given elements $a_1,
\cdots, a_m \in {\mathcal R}$ that generate ${\mathcal R}$ (as an
ideal), there exist $m \times m$ matrices $F$ and $G$ with entries
in ${\mathcal R}$ such that $a_1, \ldots ,a_m$ form the first row of
$F$ and $FG=I$ (equivalently: $GF=I$). The Hermite property of
various algebras is useful in control systems theory, see, e.g.,
\cite{Vid}.

Our central result is as follows.
\begin{thm} \label{nov182}
Let $S\subset \R^k$ be a halfspace, and let $\Sigma\subseteq S$ be
an additive semigroup: if $\lambda, \mu \in \Sigma$, then also
$\lambda+\mu\in \Sigma$. Assume $0\in\Sigma.$  Then the unital
algebras $AP^k_{\Sigma}$ and $APW^k_{\Sigma}$ are Hermite rings.
\end{thm}

The main ingredient of its proof ---
contractability of the maximal ideal spaces under the hypotheses of
Theorem \ref{nov182} --- is established in Section~\ref{mp}. Once
the contractability is proved, Theorem \ref{nov182} follows by
application of Lin's theorem \cite{Lin73}. Applications of the
result on Hermite property to rectangular factorizations and
Toeplitz corona problems are given in Sections~\ref{apfact} and
\ref{nov203}, respectively. In the preliminary Section~\ref{prel},
some known results and constructions are recalled.

Throughout the paper, if $X$ is a set (typically a Banach space or
an algebra), then we denote by $(X)^{m \times n}$ the set of $m
\times n$ matrices with entries in $X$.

\section{Auxiliary results}\label{prel}
\subsection{Hermite rings and coprime factorizations}\label{section2}
\setcounter{equation}{0}

The defining property of Hermite rings is
equivalent to its matrix analog:

\begin{lem} \label{l:1a}
If ${\mathcal R}$ is a Hermite ring, then every $k \times m$ right
(resp. $m \times k$ left) invertible matrix $F$ over ${\mathcal R}$
can be complemented, i.e., there is an $m \times m$ invertible
matrix $F_0$ over ${\mathcal R}$ such that $F$ forms the first $k$
rows (resp. $k$ columns) of $F_0$.
\end{lem}

For a proof, see p. 345 in \cite{Vid}, for example.

An important application of Hermite rings ${\mathcal R}$ without
divisors of zero has to do with coprime factorizations. Let
${\mathcal R}$ be a unital commutative ring without divisors of
zero, and let ${\mathcal F}$ be its field of fractions. Two matrices
$G_1$ and $G_2$ with entries in ${\mathcal R}$ are called {\em right
coprime} if they have same number of columns and there exist
matrices $X_1$ and $X_2$ with entries in ${\mathcal R}$ of suitable
size such that $X_1G_1+X_2G_2=I$. Dually, matrices $G_1$ and $G_2$
with entries in ${\mathcal R}$ are called {\em left coprime} if they
have same number of rows and there exist $Y_1$ and $Y_2$ with
entries in ${\mathcal R}$ such that $G_1Y_1+G_2Y_2=I$. Now let $G$
be a rectangular matrix with entries in ${\mathcal F}$. A {\em right
coprime factorization} of $G$ is, by definition, a representation of
the form $G=NM^{-1}$, where $N$ and $M$ are matrices with entries in
${\mathcal R}$, $M$ is of square size with non-zero determinant, and
$N$ and $M$ are right coprime. A {\em left coprime factorization} of
$G$ is a representation of the form $G=Q^{-1}P$, where $P$ and $Q$
are matrices with entries in ${\mathcal R}$, $Q$ is of square size
with non-zero determinant, and $P$ and $Q$ are left coprime.

In general, the existence of a coprime factorization from one side
does not imply existence of coprime factorization from the other
side. However:

\begin{prop}\label{pr:new}
Assume ${\mathcal R}$ is a Hermite ring without divisors of zero. If
$G$ is a matrix with entries in the field of fractions of ${\mathcal
R}$, and $G$ admits a left (resp., right) coprime factorization,
then $G$ admits also a right (resp., left) coprime factorization.
\end{prop}

For a proof of Proposition~\ref{pr:new} see, e.g., \cite{Vid}.

The converse of Proposition \ref{pr:new} holds as well; however, we
will not use this fact in the present paper.

\subsection{Maximal ideal spaces of commutative Banach algebras:
polynomial convexity}\label{nov171} 

 We present here some well-known background
information on commutative unital Banach algebra (all Banach
algebras are assumed to be over the complex field $\C$). The book
\cite{Stout} is used as the main reference.

Let $\CA$ be a commutative unital Banach algebra. The maximal ideal
space of $\CA$ will be denoted $\M(\CA)$. The elements of $\M(\CA)$
are nonzero multiplicative linear functionals $\phi : \CA \
\longrightarrow\ \C$. Every such functional is automatically bounded
and has norm equal to $1$. Clearly, every nonzero multiplicative
linear functional belongs to $\CA^*$, the dual space of $\CA$ of
continuous linear functionals. Thus,  $\M(\CA)\subseteq \CA^*$. The
topology in $\M(\CA)$ is induced by the weak$^*$ topology of
$\CA^*$. Note that $\M(\CA)$ is compact.

A set $\Lambda\subset \CA$ is called a {\em set of generators} of
$\CA$ if $\CA$ is the smallest closed unital subalgebra of  $\CA$
that contains  $\Lambda$. We fix a set of generators  $\Lambda$ of a
commutative unital Banach algebra  $\CA$. Let $X(\Lambda)$ be the
set of functions $f: \Lambda \ \longrightarrow \ \C$. We consider
 $X(\Lambda)$ in the product topology, i.e., in the weak$^*$ topology when
$\Lambda$ is treated as a discrete topological space.

There is a natural map
$$\chi\colon \M(\CA) \ \longrightarrow \ \mbox{$X$}(\Lambda)$$
defined by
$$ (\chi(\phi))(\lambda)=\phi(\lambda), \qquad \forall \ \ \lambda\in
\Lambda. $$
\begin{prop}
The map $\chi$ is continuous, one-to-one, and
 $\M(\CA)$ is homeomorphic to the image
$\chi( \M(\CA))$ (in the induced topology from $X(\Lambda)$).
\end{prop}
\begin{proof} The pre-base of open sets in $X(\Lambda)$
consists of the sets 
$$ \Omega_{\lambda_0,\va, f_0}:=\{f\in X(\Lambda)\, :\, |f(\zl_0)-
f_0(\zl_0)|<\va\}, \quad \zl_0\in \Lambda, \ \ f_0\in X(\Lambda), \
\ \va>0.$$ The pre-image of the set $\Omega_{\lambda_0,\va, f_0}$
under the map $\chi$ is either empty, or it consists of all elements
$\phi_0\in  \M(\CA)$ such that
$$\chi(\phi_0)\in\Omega_{\lambda_0,\va, f_0}.$$
In the latter case, together with $\phi_0$, the pre-image
$\chi^{-1}(\Omega_{\lambda_0,\va, f_0})$ contains the open set
$$ \{\phi\in \M(\CA)\, :\, |\phi(\lambda_0)-\phi_0(\lambda_0)|<
\va - |\phi_0(\lambda_0)-f_0(\lambda_0)|\}. $$ This shows that
$\chi^{-1}(\Omega_{\lambda_0,\va, f_0})$ is open. Hence $\chi$ is
continuous. If $\phi_1(\lambda)=\phi_2(\lambda)$ for every $\lambda
\in \Lambda$, where $\phi_1,\phi_2 \in  \M(\CA) $, then
$\phi_1(x)=\phi_2(x)$ for every $x\in \CA$ which can be expressed as
a finite linear combination of products of finitely many elements of
$\Lambda$. Since the set of all such linear combinations is dense in
$\CA$, by continuity of $\phi_1,\phi_2 $ we obtain that
$\phi_1(x)=\phi_2(x)$ for every $x\in \CA$. Thus, $\chi$ is
one-to-one.

Since $\M(\CA)$ is compact and $\chi$ is continuous, the image
$\chi( \M(\CA))$ is compact as well. Now clearly $\chi$ is a
homeomorphism. \end{proof}

For a fixed $\zl\in \Lambda$, denote by $\pi_{\zl}:X(\Lambda)\
\longrightarrow \ \C$ the coordinate projection on the
$\zl$-component:
$$ \pi_{\zl}(f)=f(\zl), \qquad \forall \quad f\in
X(\Lambda). $$ Clearly, $\pi_{\zl}$ is continuous.

Let $P_{\Lambda}$ be the smallest algebra (under pointwise
multiplication, addition, and scalar multiplication) of functions
$X(\Lambda)\ \longrightarrow \ \C $ that contains constant functions
and all projections $\pi_{\zl}$. Thus,  $P_{\Lambda}$ consists of
finite linear combinations of finite products of powers of the
projections $\pi_{\zl}$ (including the powers with zero exponent
that represent the constant function $1$). The functions in
$P_{\Lambda}$ are obviously continuous. If
$$ p=\sum_{j_1,\ldots, j_k\geq 0} a_{j_1,\ldots,
j_k}\pi_{\zl_{1}}^{j_1}\pi_{\zl_{2}}^{j_2} \cdots
\pi_{\zl_{k}}^{j_k},$$ where the sum is finite and $ a_{j_1,\ldots,
j_k}\in \C$,
 then for every
$f\in X(\Lambda)$ we have
\begin{equation}\label{2} p(f)=
\sum_{j_1,\ldots, j_k\geq 0} a_{j_1,\ldots,
j_k}f(\zl_{1})^{j_1}f(\zl_{2})^{j_2} \cdots f(\zl_{k})^{j_k}.
\end{equation}

The next statement is essentially Theorem 5.8 of \cite{Stout}.

\begin{prop}\label{4}
The set $\chi( \M(\CA))$ is polynomially convex, i.e.: If $f\in
X(\Lambda)$ satisfies
\begin{equation}\label{1}
 |p(f)|\leq \max\{|p(\alpha)|\, :\, \alpha\in \chi( \M(\CA))\}, \quad
\forall \quad \mbox{$p$} \in P_{\Lambda},
 \end{equation}
then $f\in \chi( \M(\CA)).$
\end{prop}

Note that because of compactness of $\chi( \M(\CA))$, the maximum is
attained in (\ref{1}). Obviously, if $f\in \chi( \M(\CA))$, then
(\ref{1}) is satisfied.

\begin{proof} Let $f\in X(\Lambda)$ be such that (\ref{1}) is
satisfied.

Define a map $h: P_{\Lambda}\ \longrightarrow \ A$ by setting
$h(\pi_{\zl})=\zl$, $h(1)=1$ (on the left here $1$ denoted the
constant function $1$), and demanding that $h$ be linear and
multiplicative. The range of $h$ is the smallest (not necessarily
closed) unital subalgebra $B$ that contains $\Lambda$. Thus $B$ is
dense in $A$.

Define $\psi \,:\, B\ \longrightarrow \ \C$ by setting
$$\psi (h(p))=p(f), \qquad \forall \quad p\in P_{\Lambda}. $$
We verify that $\psi$ is well defined. Indeed, assume $h(p)=h(q)$
for some $p,q\in P_{\Lambda}.$ Write
$$p=\sum_{j_1,\ldots, j_k\geq 0} a_{j_1,\ldots,
j_k}\pi_{\zl_{1}}^{j_1}\pi_{\zl_{2}}^{j_2} \cdots
\pi_{\zl_{k}}^{j_k},\quad q=\sum_{j_1,\ldots, j_k\geq 0}
a'_{j_1,\ldots, j_k}\pi_{\zl_{1}}^{j_1}\pi_{\zl_{2}}^{j_2} \cdots
\pi_{\zl_{k}}^{j_k}, $$ where the sums are finite, and
$$ a_{j_1,\ldots,
j_k},  a'_{j_1,\ldots, j_k} \in \C. $$

The condition $h(p)=h(q)$ means
$$ \sum_{j_1,\ldots, j_k\geq 0} a_{j_1,\ldots,
j_k}\zl_{1}^{j_1}\zl_{2}^{j_2} \cdots
\zl_{k}^{j_k}=\sum_{j_1,\ldots, j_k\geq 0} a'_{j_1,\ldots,
j_k}{\zl_{1}}^{j_1}{\zl_{2}}^{j_2} \cdots {\zl_{k}}^{j_k}. $$
Therefore, for every $\phi\in \M(\CA)$ we have $\chi(\phi)\in
X(\Lambda)$ and using formula (\ref{2}) we obtain (since $\phi$ is a
multiplicative linear functional)
$$ (p-q)(\chi(\phi))=0. $$
Now (\ref{1}) implies that
 $(p-q)(f)=0$, i.e., $p(f)=q(f)$, and $\psi$ is well defined.
Clearly, $\psi$ is a nonzero (because $\psi(1)=1$) linear
multiplicative map. If $b=h(p)\in B$, where $p\in P_{\Lambda}$, then
\begin{eqnarray*} |\psi(b)|=|p(f)|& \leq &
\max\{|p(\chi(\phi))|\, :\, \phi\in \M(\CA)\} \\ &
=&\max\{|\phi(h(p))|\, :\, \phi\in \M\CA)\} \leq
\|\mbox{$h(p)$}\|=\| \mbox{$b$}\|.
\end{eqnarray*}
Thus, $\psi$ can be extended by continuity to an element (again
denoted by $\psi$) of $\M(\CA)$. But now for every $\zl\in \Lambda$
we have
$$ f(\zl)=\pi_{\zl}(f)=\psi(h(\pi_{\zl}))=\psi(\zl), $$
and $f=\chi(\psi)$.\end{proof}

\section{Proof of the main result}\label{mp}
\setcounter{equation}{0}

We now apply the construction of Subsection~\ref{nov171} to algebras
of almost periodic functions. This will lead to the following
result:

\begin{thm} \label{nov195} Let $S\subset \R^k$ be a halfspace, and let
$\Sigma\subseteq S$ be an additive semigroup containing zero. Then
the topological spaces $\M(AP^k_{\Sigma})$ and $\M(APW^k_{\Sigma})$
are contractible.
\end{thm}

A particular case of Theorem \ref{nov195}, namely for
$\M(AP^1_{\Sigma})$, was proved in \cite{Bru05}. In our proof of the
theorem, we use some ideas and a lemma from \cite{Bru05}.

Using Lin's theorem \cite{Lin73}, \cite[Theorem 8.68]{Vid}, Theorem
\ref{nov182} is now obtained from Theorem~\ref{nov195}.

The rest of this section is devoted to the proof of
Theorem~\ref{nov195}. We focus on
the algebra $AP^k_{\Sigma}$ first.

The following property of half spaces will be needed.

\begin{prop} \label{nov172}
Let $S \subset \R^k$ be a halfspace. Then there exists a unique
vector
$Y(S)\in  \R^k$ of unit length 
such that $\langle x,Y(S)\rangle \geq 0$ for every $x\in S$.
\end{prop}

\begin{proof} For the case $S=E_k$, the standard example of a
halfspace, the proposition is obvious, with $Y(E_k)=(1,0,\ldots,
0)^T$.

We reduce the general case to the standard example. Indeed, the
halfspace $S$ induces a linear order on $\R^k$: $\lambda \preceq
\mu$, $\lambda, \mu \in \R^k$, if and only if $\mu -\lambda \in S$.
As a special case of basic results on linearly ordered vector
spaces, see \cite{Erd} or \cite[Section IV.6]{Fuchs}, we obtain that
there exists an invertible real $k \times k$ matrix $Z$ such that
$S=ZE_k:=\{A\lambda\, :\, \lambda \in E_k\}$. Now let
$Y(S)=(Z^{-1})^T Y(E_k)/\|(Z^{-1})^T Y(E_k)\|$. \end{proof}

Since for every invertible $k \times k$ matrix $Z$ the algebras
$(AP^k)_{\Sigma}$ and $(AP^k)_{Z\Sigma}$ are isometrically
isomorphic (the isometric isomorphism is induced by the map
$e_\lambda \ \ \mapsto  \ \ e_{Z\lambda}$, $\lambda\in \Sigma$, on
elementary exponentials), in view of the proof of Proposition
\ref{nov172} we may (and do) assume that $S=E_k$ and
$Y(E_k)=(1,0,\ldots, 0)^T$. We denote by $V(E_k)= \{0\}\times
\R^{k-1}$ the subspace orthogonal to $Y(E_k)$.

The set $AP^k_{\Sigma}$ of almost periodic functions $f\colon \R^k \
\longrightarrow \ \R$ with Bohr-Fourier spectrum in $\Sigma$ is a
unital commutative Banach algebra in the $L_{\infty}(\R^k)$ norm.
Furthermore,
$$\Lambda:=\{e_{\zl}\, :\zl\in \Sigma\}$$ is
obviously a generating set of $AP^k_{\Sigma}$. In the notation of
Subsection \ref{nov171}, we will work with the sets $X(\Lambda)$ of
functions $\Lambda \ \longrightarrow \ \C$ (in the weak$^*$
topology), $P_{\Lambda}$, and $$\chi(\M(AP^k_{\Sigma})) \subseteq
X(\Lambda)$$ which is homeomorphic to $\M(AP^k_{\Sigma})$.

Define the map
$$ R: X(\Lambda ) \times [0,1] \quad
\longrightarrow \quad  X(\Lambda ) $$ by
$$ R(f,t)(e_{\zl})=t^{\langle \zl, Y(E_k)\rangle}f(e_{\zl}),
\ \ \zl\in \Sigma, \quad f\in X(\Lambda ), \ \ t\in [0,1]. $$ It is
understood that here $0^0=1$.
\begin{lem} $R$ is continuous.
\end{lem}

For the proof see the proof of Lemma 2.1 in \cite{Bru05} (here the
equality $0^0=1$ is essential).

We will show that $R$ maps $\chi(\M(AP^k_{\Sigma}))\times [0,1]$
into $\chi(\M(AP^k_{\Sigma}))$.

Consider first the case when $t=0$. We have
\begin{equation}\label{3}
 R(g,0)(e_\zl)=\left\{ \begin{array}{ll} 0 & \mbox{if $\langle \zl,
Y(E_k)\rangle\neq 0$}\\[3mm]
g(e_{\zl}) & \mbox{if $\zl\in V(E_k)\cap \Sigma$} \end{array}\right.
\end{equation}
Here $g\in \chi(\M(AP^k_{\Sigma}))$. Thus, for some $\phi\in
\M(AP^k_{\Sigma})$ we have
$$ g(e_\zl)=\chi(\phi)(e_\zl)=\phi(e_\zl), \quad \zl\in \Sigma. $$

Let $APP^k_\Sigma$ be the linear set of all almost periodic
polynomials (i.e., almost periodic functions whose Bohr-Fourier
spectrum is a finite set) of several variables with Bohr-Fourier
spectrum in $\Sigma$. Define the map
$$ \phi'\colon APP^k_{\Sigma} \ \longrightarrow \ \C $$
by
$$ \phi'(e_\zl)=\left\{ \begin{array}{ll} \phi (e_\zl)
& \mbox{if $\zl \in V(E_k)\cap \Sigma$,} \\[3mm]  0 &
\mbox{if $\zl \in \Sigma\setminus V(E_k)$}, \end{array}\right. $$
and extend it by linearity to $(APP^k)_{\Sigma} $. The map $\phi'$
is clearly unital. It is also multiplicative, because $\phi$ is
multiplicative and in view of the property
$$ \zl_1+\zl_2\in  \Sigma\setminus V(E_k) $$
provided $\zl_1,\zl_2\in \Sigma$ and at least one of $\zl_1,\zl_2$
does not belong to $V(E_k)$. Moreover, $\phi'$ is bounded (in the
$\|\cdot\|_\infty$ norm). Indeed, $\phi'$ is a composition of the
bounded linear functional $\phi$ and the linear projection $P$
defined on $APP^k_\Sigma$ by the equality
$$ P(e_\lambda)=\left\{ \begin{array}{ll} e_\zl
& \mbox{if $\zl \in V(E_k)\cap \Sigma$,} \\[3mm]  0 &
\mbox{if $\zl \in \Sigma\setminus V(E_k)$}, \end{array}\right.
$$
It is easy to see that
$$ P(f)=\lim_{T \ \longrightarrow  \ \infty}
\ \frac{1}{2T} \, \int_{-T}^T f(x,x_2,\ldots, x_k)dx, \quad f\in
APP^k_\Sigma. $$ Thus, $P$ is bounded, and hence so is $\phi'$.
Therefore, $\phi'$ extends by continuity to an element (also denoted
by  $\phi'$)  of $\M(AP^k_{\Sigma})$. Formula (\ref{3}) shows that
$$ R(g,0)=\chi(\phi'),$$
and therefore
$$ R(g,0)\in  \chi(\M(AP^k_{\Sigma})). $$

Consider now the case $t>0$. Arguing by contradiction assume that
there is a $v\in \chi(\M(AP^k_{\Sigma}))$ and $t>0$ such that
$R(v,t)\not\in \chi(\M(AP^k_{\Sigma})).$ By Proposition~\ref{4},
there is a polynomial
$$ p=\sum_{j_1,\ldots, j_k\geq 0} a_{j_1,\ldots,
j_k}\pi_{\zl_{1}}^{j_1}\pi_{\zl_{2}}^{j_2} \cdots
\pi_{\zl_{k}}^{j_k}\in P_{\Lambda}, \quad \zl_1,\ldots, \zl_k\in
\Sigma,$$ such that
$$ |p(R(v,t))|> \max
\{|p(\alpha)|\, :\, \alpha\in \chi( \M((AP^k)_{\Sigma})\}. $$ We
have
$$ p(R(v,t))=
\sum_{j_1,\ldots, j_k\geq 0} a_{j_1,\ldots, j_k}t^{\langle
j_1\zl_1+\cdots +j_k\zl_k,
Y(E_k)\rangle}v(e_{\zl_{1}})^{j_1}v(e_{\zl_{2}})^{j_2} \cdots
v(e_{\zl_{k}})^{j_k}. $$ Define the polynomial
$$ q:=\sum_{j_1,\ldots, j_k\geq 0} a_{j_1,\ldots,
j_k}t^{\langle j_1\zl_1+\cdots +j_k\zl_k,
Y(E_k)\rangle}\pi_{\zl_{1}}^{j_1}\pi_{\zl_{2}}^{j_2} \cdots
\pi_{\zl_{k}}^{j_k}\in P_{\Lambda}.$$ Then $ p(R(v,t))=q(v)$, and we
have
\begin{equation}\label{oct306} |q(v)| >
\max \{|p(\alpha)|\, :\, \alpha\in \chi( \M(AP^k_{\Sigma}))\}.
\end{equation}
Consider for a fixed $t$, $0<t\leq 1$, the following transformation:
If $\alpha=\chi(\phi)$ for some $\phi\in \M(AP^k_{\Sigma})$, we let
$\alpha_t=\chi(\phi_t)$, where
$$\phi_t\colon (APP^k)_{\Sigma} \ \longrightarrow \ \C$$ is defined by the
property that
$$\phi_t(e_{\zl})=t^{\langle\zl, Y(E_k)\rangle}\phi(e_{\zl}), \quad \forall \quad
\zl\in \Lambda, $$ and extended by linearity to the set of almost
periodic polynomials in $AP^k_{\Sigma}$. The map is also
multiplicative, because
$$ \phi_t(e_{\zl}e_{\mu})=\phi_t(e_{\zl+\mu})=t^{\langle\zl+\mu,
Y(E_k)\rangle}\phi(e_{\zl+\mu}) =\phi_t(e_{\zl})\phi_t(e_{\mu})$$
for all $\zl,\mu\in \Sigma.$

Next, we show that $\phi_t$ is bounded. To this end, we will use a
proposition:
\begin{prop}\label{mar211}
Let $f\in AP^k_{\Omega}$, where
$$\Omega \, \subseteq \, {\rm closure}\ {\rm of} \  E_k
:=\{(\lambda_1,\ldots, \lambda_k)\in \R^k\, :\, \lambda_1\geq 0\}.
$$ Then, for every fixed $(k-1)$-tuple $(x'_2, x'_3, \ldots,
x'_k)\in \R^{k-1}$, the function $\widetilde{f}$ defined by
$\widetilde{f}(x'):=f(x',x'_2,\ldots, x'_k)$ has the following
properties:
\begin{itemize}
\item[${\rm (a)}$] $\widetilde{f}$
is almost periodic with Bohr-Fourier spectrum in $\R_+$, the set of
nonnegative real numbers;
\item[${\rm (b)}$] $\widetilde{f}$ admits a unique continuation, also denoted by
$\widetilde{f}$, into the open complex upper half-plane $\C_+$ such
that the continuation of $\widetilde{f}$ is continuous on $\C_+\cup
\R_+$ and analytic on $\C_+$;
\item[${\rm (c)}$]  For every fixed $y'>0$, the
slice $\widetilde{f}_{y'}(x'):= \widetilde{f}(x'+iy')$ is an almost
periodic function of $x'$, and
 $$ \|\widetilde{f}_{y'}\|_\infty\leq \|\widetilde{f}\|_\infty. $$
\item[${\rm (d)}$] For every fixed $y'>0$, the function $\widehat{f}_{y'}(x',x'_2,\ldots, x'_k):=
\widetilde{f}(x'+iy')$ is an almost periodic function of the
variables $(x', x'_2,\ldots, x'_k)$, and $\sigma
(\widehat{f}_{y'})\subseteq \Omega$.
\item[${\rm (e)}$] If in addition $f\in APW^k_{\Omega}$, then
for every fixed $y'>0$, we have:
\begin{itemize}
\item[]
$\widetilde{f}_{y'}\in APW^1_\Omega$;
\item[] the inequality
$\|\widetilde{f}_{y'}\|_W\leq \|\widetilde{f}\|_W $ holds;
\item[] The function  $\widehat{f}_{y'}$ belongs to
$APW^k_\Omega$ as a function of $(x',x'_2,\ldots, x'_k)\in \R^k$.
\end{itemize} \end{itemize}\end{prop}

\begin{proof} Part (a) follows easily from the fact that, given any
set $\Delta\subseteq \R^k$, every function in $AP^k_{\Delta}$ can be
uniformly approximated by almost periodic polynomials with
Bohr-Fourier spectrum  in $\Delta$; this goes back to \cite{Bes}.
Part (b) is standard given that $\sigma(\widetilde{f})\subset \R_+$
by (a). For (d), we use the Poisson formula (see \cite[Chapter
8]{Hoffman} or \cite[Chapter 5]{RosRov97}, for example):
\begin{equation}\label{mar212} \widetilde{f}_{y'}(x')=\frac{y'}{\pi}\int_{-\infty}^{\infty}
\frac{f(t,x'_2,\ldots, x'_k)}{(t-x')^2+(y')^2}dt, \quad y'>0.
\end{equation} Let $h^{(1)},h^{(2)},\ldots $ be a sequence of almost
periodic polynomials with $\sigma (h^{(j)})\subseteq \Omega$ and
such that
$$ \sup_{x\in \R^k} |h^{(j)}(x)-f(x)|<\frac{1}{j}, \qquad j=1,2,\ldots .$$
For every fixed $y'>0$, the functions
$\widehat{h^{(j)}}_{y'}(x',x'_2,\ldots, x'_k)$, $j=1,2,\ldots$, are
obviously almost periodic polynomials with Bohr-Fourier spectrum in
$\Omega$. Now (\ref{mar212}) gives for a fixed $x'\in \R$:
$$ \sup_{(x'_2,\ldots, x'_k)\in \R^{k-1}}
\, |\widehat{f}_{y'}(x',x'_2,\ldots, x'_k)-
\widehat{h^{(j)}}_{y'}(x',x'_2,\ldots, x'_k)|\leq
$$ \begin{equation}\label{mar213} \frac{y'}{\pi}\int_{-\infty}^{\infty}
\frac{\sup_{(x'_2,\ldots, x'_k)\in \R^{k-1}}\,|f(t,x'_2,\ldots,
x'_k)-h^{(j)}(t,x'_2,\ldots, x'_k)|}{(t-x')^2+(y')^2}dt \ \leq \
\frac{1}{j},\end{equation} hence
$$ \sup_{(x',x'_2,\ldots, x'_k)\in \R^{k}}
\, |\widehat{f}_{y'}(x',x'_2,\ldots, x'_k)-
\widehat{h^{(j)}}_{y'}(x',x'_2,\ldots, x'_k)| \leq \frac{1}{j}, $$
and (d) follows. Part (c) can be proved by using a string of
inequalities analogous to (\ref{mar213}). Finally, the part (e) is
easy to verify taking into account that every almost periodic
function in the Wiener algebra is the sum of its absolutely
convergent Bohr-Fourier series. \end{proof}

We have
$$\phi_t(e_{\zl})=\phi(t^{\langle\zl, Y(E_k)\rangle}e^{i\langle\zl,
x\rangle})=\phi(e^{i\langle\zl,x-i\log t Y(E_k)\rangle}). $$ Note
that $-\log t \geq 0$ because $t\leq 1$. So for any almost periodic
polynomial $f$ in $AP^k_{\Sigma}$ we have
$$ \phi_t(f(x))=\phi(f(x-i\log t Y(E_k))). $$
(Note that the right hand side is well defined because of
Proposition \ref{mar211}(d).) But
$$ \sup_{x\in \R^k} |f(x-i\log t Y(E_k))|
=\sup_{(x_2,\ldots, x_k)\in \R^{k-1}} \ \left\{ \sup_{x'\in \R}
|f(x'- i\log t,x_2, \ldots ,x_k)|\right\} $$ $$ \leq
\sup_{(x_2,\ldots, x_k)\in \R^{k-1}} \ \left\{ \sup_{x'\in \R}
|f(x',x_2, \ldots ,x_k)|\right\}= \sup_{x\in \R^k} |f(x)|, $$ where
the inequality follows from Proposition \ref{mar211}(c). This proves
that $\phi_t$ is bounded on the set of almost periodic polynomials.
Now we can extend $\phi_t$ by continuity to a bona fide $\phi_t\in
 \M(AP^k_{\Sigma})$. From (\ref{oct306}) we have
\begin{eqnarray*}  |q(v)| & > &
\max \{|p(\alpha_t)|\, :\, \alpha\in \chi( \M(AP^k_{\Sigma}))\}\\
&= & \max\{|q(\alpha)|\, :\, \alpha\in \chi( \M(AP^k_{\Sigma}))\},
\end{eqnarray*}
a contradiction with $v\in \chi( \M(AP^k_{\Sigma}))$.

We have proved that $R$ maps $\chi(\M(AP^k_{\Sigma}))\times [0,1]$
into $\chi(\M(AP^k_{\Sigma}))$. Note that  $R$ is continuous, when
restricted to $\chi(\M(AP^k_{\Sigma}))\times \{1\}$ the map $R$ is
the identity, and the range of its restriction to
$\chi(\M(AP^k_{\Sigma}))\times \{0\}$ is the maximal ideal space of
the algebra $AP^k_{\Sigma\cap V(E_k)}$. The latter algebra may be
identified with $AP^k_{\Sigma'}$, where $\Sigma'$ is a suitable
semigroup of $\R^{k-1}$; in other words, we identify $V(E_k)$ with
$\R^{k-1}$, and note that $E_k \cap V(E_k)$ is a halfspace in
$V(E_k)$. Now we apply an analogous procedure in the case of $k-1$
variables. After $k$ such steps, we end up with a continuous map
$$\widetilde{R}\colon \chi(\M(AP^k_{\Sigma}))\times [0,1] \ \longrightarrow
\ \chi(\M(AP^k_{\Sigma})) $$ such that $\widetilde{R}$ is
continuous, the map $\widetilde{R}$ is the identity when restricted
to $\chi(\M(AP^k_{\Sigma}))\times \{1\}$, and the range of its
restriction to $\chi(\M(AP^k_{\Sigma}))\times \{0\}$ is a singleton.
This proves Theorem \ref{nov195} for the algebra $AP^k_{\Sigma}$.

For the contractibility of $\M(APW^k_{\Sigma})$ the proof is
analogous; we have to use part (e) of Proposition \ref{mar211} in
that case.

\section{AP Factorizations of Rectangular
Matrix Functions}\label{apfact} \setcounter{equation}{0}

In this section we apply the Hermite properties established in
Theorem \ref{nov182} to factorizations of almost periodic matrix
functions of several variables.

We start with the concept of factorization. Throughout this section,
we fix a halfspace $S\subset \R^k$ and an additive subgroup $\Sigma
\subseteq \R^k$. Let $G\in (AP^k_\Sigma)^{m\times n}$. A {\em left}
$(AP_S)_\Sigma$ {\em factorization} of $G$ is a representation of
the form
\begin{equation}\label{factor}
 G=G^+{\diag} (e_{\lambda_1}, \ldots ,e_{\lambda_p}) G^-, \end{equation}
where $G^+\in (AP^k_{\Sigma\cap S})^{m \times p}$, $G^-\in
(AP^k_{\Sigma\cap (-S)}) ^{p \times n}$, $G^+$ has a right inverse
in $(AP^k_{\Sigma\cap S})^{p \times m}$, $G^-$ has a left inverse in
$(AP^k_{\Sigma\cap (-S)}) ^{n \times p}$, and $\lambda_1,\ldots,
\lambda_p \in \Sigma$. The elements $\lambda_j$'s are called the
{\em factorization indices}; they are unique if we require
$\lambda_j-\lambda_{j-1}\in S$ for $j=2,\ldots,n$ (we will assume
that these inclusions hold). If $APW$ is used in place of $AP$
throughout in the above definition, then we say that (\ref{factor})
is a  {\em left} $(APW_S)_\Sigma$ {\em factorization} of a matrix
function $G\in  (AP^k_\Sigma)^{m\times n}$. We say that a left
$(AP_S)_{\Sigma}$ (or $(APW_S)_{\Sigma}$) factorization
(\ref{factor}) is {\it canonical} if the factorization indices are
zero, that is, if the middle factor in (\ref{factor}) is the
identity matrix $I_p$. A right  $(AP_S)_{\Sigma}$ (or
$(APW_S)_{\Sigma}$) factorization is introduced analogously, with
the roles of $G_+$ and $G_+$ interchanged, in other words, a right
 $(AP_S)_{\Sigma}$ (or $(APW_S)_{\Sigma}$) factorization coincides with
a left  $(AP_{-S})_{\Sigma}$ (or $(APW_{-S})_{\Sigma}$)
factorization. Therefore, we focus on the left factorization only,
and the adjective ``left" will be often omitted. A matrix function
$G\in (AP^k_\Sigma)^{m\times n}$ is said to be
$(AP_S)_{\Sigma}$-{\em factorable} if it admits an
$(AP_S)_{\Sigma}$-factorization; analogously
$(APW_S)_{\Sigma}$-factorable matrix functions are introduced.

Left invertibility of $G^-$ and right invertibility of $G^+$ in the
formula (\ref{factor}) imply that $p\geq\max\{ m,n\}$. On the other
hand, if no upper bound is imposed on $p$ then every matrix function
$G\in (APW^k_\Sigma)^{m\times n}$ can be represented as $G=G^+G^-$
with
\begin{equation}\label{triv}
G^+= \left[ \Pi_{S}G \ \ \ I_m\right], \quad G^- =
\left[\begin{matrix} I_n\\ \Pi_{S\setminus\{0\}}G\end{matrix}\right]
,
\end{equation}
where for a nonempty set $\Lambda \subseteq \R^k$ we denote by
$\Pi_\Lambda$ the projection
\begin{equation} \label{prl}
\Pi_\Lambda f = \sum_{\lambda\in\Lambda\cap\sigma(f)}f_\lambda
e_\lambda, \qquad f\in (APW^k)^{m \times n}.
\end{equation}
To exclude the triviality (\ref{triv}), we impose an additional
condition
\begin{equation}\label{pmn}
p = \max\{ m,n\}
\end{equation}
on the factorization (\ref{factor}).

Under condition (\ref{pmn}), for square ($m=n$) matrix functions $G$
the factors $G^\pm$ in (\ref{factor}) are also square and therefore
two sided invertible in the algebras $(AP^k_{\Sigma\cap (\pm
S)})^{n\times n}$.

If $k=1$, the only half spaces are $\R_\pm :=\{ x\in \R\colon \pm
x\geq 0\}$. Therefore, representation (\ref{factor}) with $m=n=p$
and $\Sigma=\Z$ up to a simple change of variable $\zeta=\exp (2\pi
ix)$ is the classical Wiener-Hopf factorization on the unit circle
$\T$. The pioneer work on this kind of factorization was done by
I.~Gohberg and M.~Krein \cite{GK58}, see e.g. \cite{GF74,CG,GKS03}
for the further development. The Wiener-Hopf factorization of
rectangular matrix functions with $G_+$ ($G_-$) being left
(respectively, right) invertible was treated in \cite{RakSpit96}.

The case $k=1$, $\Sigma=\R$ (again, with $m=n$) corresponds to the
$AP$ (or $APW$) factorization of matrix functions, introduced by
Yu.~Karlovich and one of the authors (see e.g. \cite{KarlSpit89})
and considered further in subsequent publications
\cite{KarlSpit951,KarlSpit952,BKST,QRS,BKdST1,BKdST2}. A systematic
exposition of the $AP$ factorization and its properties can be found
in \cite{BKS1}. $AP$ factorizations with respect to arbitrary
subgroups $\Sigma$ of $\R$ were studied in \cite{RSW98}, and the
transition to $k>1$ was accomplished in \cite{RSW01ot,RSW03}. These
papers, along with \cite{BKRS,RSW02jot,RSWmem}, contain also
applications of the $(AP_S)_\Sigma$ factorization to various
extension and interpolation problems, as well as spectral
estimation. We mention here for completeness that $(APW_S)_\Sigma$
factorization can be thought of as a particular case of the
factorization in Wiener algebras in the setting of ordered Abelian
groups; see \cite{MRSW,MRS05,RS06} and \cite{EMRS} for recent
developments in this direction.

To formulate a factorization criterion for square matrix functions,
let us recall that the formula
\begin{equation}\label{inner}
\scal{f,g} = M\{fg^*\}, \ f,g \in (AP^k)
\end{equation}
defines an inner product on $AP^k$. The completion of $AP^k$ with
respect to this inner product is called the {\em Besicovitch space}
and is denoted by $B^k$. Thus $B^k$ is a Hilbert space. The
projection (\ref{prl}) extends by continuity to the orthogonal
projection (also denoted $\Pi_\Lambda$) on $B^k$. We denote by
$B^k_\Lambda$ the range of $\Pi_\Lambda$, or, equivalently, the
completion of $AP^k_\Lambda$ with respect to the inner product
(\ref{inner}). The vector valued Besicovitch space $(B^k)^{n \times
1}$ consists of $n \times 1$ columns with components in $B^k$, with
the standard Hilbert space structure. Similarly, $(B^k_\Lambda)^{n
\times 1}$ is the Hilbert space of $n \times 1$ columns with
components in $B^k_\Lambda$.

For $F\in (AP^k_{\Sigma})^{m\times n}$, where $\Sigma$ is an
additive subgroup of $\R^k$, the {\em Toeplitz operator}
$T(F)_{S,\Sigma}$ acts from $\left( B^k_{S\cap\Sigma}\right)^{n
\times 1}$ to $\left( B^k_{S_\cap\Sigma}\right)^{m \times 1}$
according to the formula
\[
T(F)_{S,\Sigma}\phi = \Pi_{S}(F\phi) =\Pi_{S\cap \Sigma}(F\phi).
\]
\begin{thm} \label{th:inv}
Let $\Sigma$ be an additive subgroup of $\R^k$ and let $G\in
(APW^k_\Sigma)^{n\times n}$. Then $G$ admits a canonical
$(AP_S)_{\R^k}$ factorization if and only if the operator
$T(G^T)_{S,\R^k}$ is invertible, where the superscript $T$ indicates
the transposed matrix. If this is the case, then the operator
$T(G^T)_{S,\Sigma}$ is invertible as well, and any $(AP_S)_{\R^k}$
factorization of $G$ automatically is its canonical $(APW_S)_\Sigma$
factorization.
\end{thm}

Theorem~\ref{th:inv} is a subset of a lengthy set of equivalent
statements constituting Theorem~2.3 of \cite{BKRS} and based on
previous results from \cite{Spit89,Karl90,Karl93,RSW98,BKdST2}. We
refer an interested reader to \cite{BKRS} for the complete proof and
the detailed history of the matter. 

For non-square matrix functions, the following result holds.
\begin{thm} \label{th:recfac}
Let $G\in (AP^k)_\Sigma^{m\times n}$ with $m<n$ (respectively,
$m>n$). Then $G$ is $(AP_S)_\Sigma$-factorable with indices
$\lambda_1,\ldots, \lambda_p$ if and only if $G$ can be augmented by
$n-m$ rows (respectively, $m-n$ columns) to a square
$(AP_S)_\Sigma$-factorable matrix function with the same indices. In
particular, $G$ admits a canonical $(AP_S)_\Sigma$-factorization if
and only if $G$ can be augmented to a square size matrix function
that admits a canonical $(AP_S)_\Sigma$-factorization.
\end{thm}
\begin{proof} For the sake of definiteness let $m<n$; the case
$m>n$ can be considered similarly.

Sufficiency: Suppose that $F$ is a square $n\times n$ matrix the
first $m$ rows of which form $G$: $F=\left[\begin{matrix} G \\
F_0\end{matrix}\right]$. Let $F=F^+DF^-$ be an
$(AP_S)_\Sigma$-factorization of $F$. Then $G=G^+DF^-$, where $G^+$
is formed by the first $m$ rows of $F^+$, is an
$(AP_S)_\Sigma$-factorization of $G$.

Necessity: Let $G=G^+DG^-$ be an $(AP_S)_\Sigma$-factorization of
$G$. According to Theorem \ref{nov182} and Lemma \ref{l:1a}, a right
invertible matrix function $G^+$ can be row-augmented to an
invertible in $(AP^k)_{S\cap\Sigma}^{n \times n}$ square matrix
$F^+$. Let $F=F^+DG^-$. Then $F$ is a row augmentation of $G$ which
is obviously $(AP_S)_\Sigma$-factorable. \end{proof}

An $(APW^k_\Sigma)^{m \times n}$ analog of Theorem~\ref{th:recfac}
is also valid, with essentially the same proof. Combining this
analog (for the special case of canonical factorization) with
Theorem \ref{th:inv}, we obtain:
\begin{cor}
The following statements are equivalent for $G\in
(APW^k_\Sigma)^{m\times n}$:
\begin{itemize}
\item[${\rm (1)}$] $G$ admits a canonical $(APW_S)_{\R^k}$
factorization;
\item[${\rm (2)}$]
$G$ can be augmented to a square size matrix function that admits a
canonical $(APW_S)_{\R^k}$-factorization;
\item[${\rm (3)}$]  $G$ can be augmented to a matrix function
$\widetilde{G}\in (APW^k_\Sigma)^{\max\{m,n\}\times \max\{m,n\}}$
such that the operator $T(\widetilde{G}^T)_{S,\R^k}$ is invertible.
\end{itemize}
In each of the statements ${\rm (1) - (3)}$, $\R^k$ can be replaced
with $\Sigma$.
\end{cor}

For an $(APW_S)_\Sigma$-factorization to exist it is obviously
necessary that
\begin{equation}\label{mar151} G\in (APW^k)_\Sigma^{m\times n} \quad \mbox{and}\quad
\mbox{$G$ has one-sided inverse in $(APW^k)\Sigma^{n\times m}$}.
\end{equation} These conditions are sufficient for existence of an
$(APW_S)_\Sigma$-factorization in the case $m=n=1$, but it is well
known that for $m=n>1$ these conditions are not sufficient, see
\cite{KarlSpit89} for the first example of this kind. We do not know
what the situation is with regard to sufficiency of (\ref{mar151})
when $m\neq n$, in particular, when $\min\{ m,n\}=1<\max\{m,n \}$.
However, if in the latter case one of the entries of $G$ is
invertible then $G$ is $(APW_S)_\Sigma$-factorable. If, for example,
$G=\left[ g_1,\ldots , g_n\right]$ with $g_j\in (APW^k)_\Sigma$,
$j=1,2,\ldots, n$, and invertible $g_1$, then
\[
G=\left[ g_1^+, h_2^+,\ldots, h_n^+\right] \ (e_\lambda I_n)
\left[\begin{matrix} g_1^- & h_2^- & \ldots & h_n^- \\ & \ddots & & \\
& & & \\
0 & \ldots & 0 & g_1^-\end{matrix}\right]
\]
is the desired factorization. Here $g_1 = g_1^+ e_\lambda g_1^-$ is
an $(AP_S)_\Sigma$-factorization of the scalar function $g_1$; note
that an $(AP_S)_\Sigma$-factorization of a scalar almost periodic
function with Bohr-Fourier spectrum in $\Sigma$ exists provided the
function is invertible, i.e., takes nonzero values bounded away from
zero, a result proved in \cite{EMRS} (see \cite{EhrMee03} for
special cases). Furthermore,
$$h_j^+=g_1^+\left(\Pi_{S} g_1^{-1}\right)g_j, \quad j=2,3,\ldots, n, $$ and
$$h_j^-=g_1^-\left(\Pi_{(-S)\setminus \{0\}} g_1^{-1}\right)g_j, \quad j=2,\ldots,n. $$

\section{The Toeplitz Corona Problem}\label{nov203}
\setcounter{equation}{0}

To state the main result of this section, we fix a halfspace
$S\subset \R^k$ and denote by $\|X\|$ the operator norm of an $p
\times q$ matrix $X$: $\|X\|=\max_{x\neq 0} \left\{\|Xx\|_2/
\|x\|_2\right\}$, where $\|x\|_2$ is the standard Euclidean norm in
$\C^q$. Also, we will use the notion of left coprimeness given in
Subsection~\ref{section2}.

\begin{thm} \label{T:genToepcor}
Suppose that we are given $A \in (APW^k_S)^{p \times m}$, $B \in
(APW^k_S)^{p \times p}$ with $B$ invertible in $(APW^k)^{p \times
p}$, and a positive number $\gamma$.  Assume in addition that $A$
and $B$ are left coprime over $APW^k_{S}$.   Let $\Lambda'$ be any
additive subgroup of ${\mathbb R}^k$ containing $\sigma(A) \cup
\sigma(B)$.  Under these conditions there exists  an $F \in
(APW^k_{S \cap \Lambda^{\prime}})^{m \times p}$ such that
\begin{equation}\label{solution}
\|F\|_{\infty}:=\sup_{t\in \R^k} \|F(t)\| \le \gamma, \quad {\rm
and}\quad AF = B
\end{equation}
only if
\begin{equation} \label{Toepineq}
 T(A)^{}_{S,\Lambda^{\prime}}  (T(A)^{}_{S, \Lambda^{\prime} })^{*}
 \ge \frac{1}{\gamma^{2}} T(B)^{}_{S,\Lambda^{\prime}}
 (T(B)^{}_{S,\Lambda^{\prime}})^{*}.
\end{equation}

Conversely, if the operator
\begin{equation} \label{Toepnew}
T(A)^{}_{S,\Lambda^{\prime}}  (T(A)^{}_{S, \Lambda^{\prime} })^{*} -
\frac{1}{\gamma^{2}} T(B)^{}_{S,\Lambda^{\prime}}
 (T(B)^{}_{S,\Lambda^{\prime}})^{*}
\end{equation}
is positive definite and invertible, then there exists $F\in
(APW^k_{S\cap \Lambda^{\prime}})^{m \times p}$ such that {\em
(\ref{solution})} holds. The family of all such matrix functions $F$
can be parameterized by $G\in \left(APW^k_{\Lambda^{\prime} \cap
S}\right)^{(m-p) \times p}$  with $\| G \|_{\infty} \le 1$ in terms
of a certain fractional linear transformation
\begin{equation}\label{apr31}  F =
(\Theta_{11}G + \Theta_{12})(\Theta_{21} G + \Theta_{22})^{-1}.
\end{equation}
\end{thm}

The proof of Theorem \ref{T:genToepcor} gives a formula (see formula
(\ref{apr32}))
for the function $\Theta=\left[\begin{array}{cc} \Theta_{11} & \Theta_{12}\\
\Theta_{21} & \Theta_{22}\end{array}\right]$ that appears in
(\ref{apr31}).
\medskip

In connection with Theorem~\ref{T:genToepcor} notice that $B\in
(APW^k)^{p \times p}$ is invertible in $(APW^k)^{p \times p}$ if and
only if $B$ is invertible in $(APW^k_{\Lambda})^{p \times p}$, where
$\Lambda$ is the additive subgroup of $\R^k $ generated by $\sigma
(B)$, if and only if $|\det (B(t))|\geq \epsilon >0$ for all $t\in
\R^k$, where $\epsilon $ is independent of $t$. See, e.g.,
Proposition 2.2 and Corollary 2.7 of \cite{RSW98} for a proof of
this statement.

The problem (considered in Theorem~\ref{T:genToepcor}) of solving
the functional equation $AF=B$ for given functions $A$ and $B$ under
the additional restriction $\|F\|\leq \gamma$ for some $\gamma>0$ is
known as the {\em Toeplitz corona problem}. It has been extensively
studied in the literature, in various settings. See, for example,
\cite{Hel86,Ros} for the background on the Toeplitz corona problem
in the context of the algebra $H^\infty$ of the upper half plane or
of the unit disk. Results on the Toeplitz corona problem in
polydisks, and more generally, in reproducing kernel Hilbert space
are found in \cite{BaTr,btv}. For algebras of almost periodic
functions, the special case of Theorem ~\ref{T:genToepcor}, where $B
= I$ and $k=1$, appears as part of Theorem 5.2 in \cite{BKRS}.
Theorem~4.1 of \cite{BallRS} is another version of Theorem
\ref{T:genToepcor}, also treating the case $k=1$.

The proof of the general case of Theorem \ref{T:genToepcor} follows
using an approach similar to that of \cite{BKRS,BallRS}, once we
verify the following lemma.

\begin{lem} \label{L:kerran}  Let $\Lambda'$ be a subgroup of $\R^k$.
Suppose that the pair of matrix functions $(A,B)$, with $A \in
(APW^k_{\Lambda^{\prime}\cap S})^{p \times m}$, $B \in
(APW^k_{\Lambda^{\prime}\cap S})^{p \times p} $ and $B$ invertible
in $\left(APW^k_{\Lambda^{\prime}}\right)^{p \times p}$, is a left
coprime pair over $(APW^k)_{\Lambda^{\prime} \cap S}$. Then:

${\rm (i)}$  the function $W = B^{-1}A$ has a right coprime
factorization $W = C D^{-1}$ over $(APW^k)_{\Lambda^{\prime} \cap
S}$, with $$C \in  (APW^k)^{p \times m}_{\Lambda^{\prime} \cap S},
\qquad D \in (APW^k)^{m \times m}_{\Lambda^{\prime} \cap S}, $$ and
$D$ invertible in $(APW^k_{\Lambda^{\prime}})^{m \times m}$.

Moreover,

${\rm (ii)}$ the following  kernel-range Toeplitz operator identity
holds:
\begin{equation}  \label{kerran}
 {\rm Ker }\,  \begin{bmatrix} T(A)^{}_{S,\Lambda^{\prime}} &
 -T(B)^{}_{S,\Lambda^{\prime}} \end{bmatrix}
= {\rm Im }\, \begin{bmatrix} T(D)^{}_{S,\Lambda^{\prime}} \\[2mm]
T(C)^{}_{S,\Lambda^{\prime}}
\end{bmatrix}.
\end{equation}
\end{lem}

\begin{proof} By assumption,  the pair $(A,B)$ is left coprime, and
hence $W = B^{-1}A$ is a left coprime factorization. It follows from
Theorem \ref{nov182} that $APW^k_{S\cap \Lambda^{\prime}}$ is a
Hermite ring. It has no divisors of zero, as easily follows from the
analytic continuation property of Proposition \ref{mar211}. By
Proposition~\ref{pr:new} $W$ admits a right coprime factorization.
Hence we may write $B^{-1}A = W = C D^{-1}$ where $C \in
(APW^k_{\Lambda^{\prime} \cap S})^{p \times m}$ and $D \in
(APW^k_{\Lambda^{\prime} \cap S})^{m \times m}$ with the determinant
of $D$ not identically zero. The rest of the proof goes exactly as
that of \cite[Lemma 4.2]{BallRS}. \end{proof}

We also need the following property of canonical factorizations.

\begin{thm}\label{nov206}
If a Hermitian matrix function admits a canonical factorization,
then a factorization can be obtained in the symmetric form
$A=B^*DB$, where $D$ is a diagonal matrix with $\pm 1$'s on the main
diagonal.
\end{thm}

The classical version of this result (for the Wiener-Hopf
factorization) goes back to \cite{Shmu54}, the $AP$ version (for
$k=1$) is in \cite{Spit89}, see also \cite[Corollary 9.13]{BKS1}.
Finally, the general case is in \cite[Theorem 5.1]{RSW01ot}.

{\sl Proof of Theorem \ref{T:genToepcor}.} Necessity of the
condition \eqref{Toepineq} is straightforward. Indeed, if $F \in
(APW^k_{\Lambda^{\prime} \cap S })^{m \times p}$ satisfies $\|
F\|_{\infty} \le \gamma$ and $A F = B$, then
 $$ \begin{array}{rcl}
 T(B)^{}_{S,\Lambda^{\prime}}
 (T(B)^{}_{S,\Lambda^{\prime}})^{*} & = &
  T(A)^{}_{S,\Lambda^{\prime}}
   T(F)^{}_{S,\Lambda^{\prime}}
( T(F)^{}_{S,\Lambda^{\prime}})^{*}
 (T(A)^{}_{S,\Lambda^{\prime}})^{*}  \\[3mm]
 & \le & \gamma^{2} T(A)^{}_{S,\Lambda^{\prime}}
 (T(A)^{}_{S,\Lambda^{\prime}})^{*}
\end{array}
$$
and \eqref{Toepineq} follows.

We next consider the proof of the converse statement.  First observe
that, by replacing $B$ with $\gamma^{-1}B$ and $F$ with
$\gamma^{-1}F$, we may assume without loss of generality that
$\gamma = 1$. The assumption that the operator
$$T(A)^{}_{S,\Lambda^{\prime}} ( T(A)^{}_{S,\Lambda^{\prime}})^{*}
 -  T(B)^{}_{S,\Lambda^{\prime}}
 (T(B)^{}_{S,\Lambda^{\prime}})^{*}$$ is positive definite
and invertible has the geometric interpretation that the subspace
$$
{\mathcal P} := \text{Im } \begin{bmatrix}  (T(A)^{}_{S,\Lambda^{\prime}})^{*} \\[2mm]
(T(B)^{}_{S,\Lambda^{\prime}})^{*} \end{bmatrix}
$$ is a strictly positive subspace in the $J_{1}$-inner product on
$\left(B^k_{\Lambda^{\prime}\cap S}\right)^{(m+p) \times 1}$, where
$$J_{1} =
\begin{bmatrix} I_{m} & 0 \\ 0 & -I_{p} \end{bmatrix} .$$
(Here and in the sequel we use well known basic properties of
geometry of Krein spaces, see, for example, \cite{AI}.) Consequently
the $J_{1}$-orthogonal complement ${\mathcal P}^{\perp J_{1}}$ of
${\mathcal P}$ is a regular subspace of
$\left(B^k_{\Lambda^{\prime}\cap S}\right)^{(m +p) \times 1}$ in the
Krein space with the $J_1$-inner product. One can easily verify that
\begin{equation}\label{pj1}
{\mathcal P}^{\perp J_{1}} = \text{Ker }
 \begin{bmatrix}  T(A)^{}_{S,\Lambda^{\prime}} &
 - T(B)^{}_{S,\Lambda^{\prime}} \end{bmatrix}.
 \end{equation}
Indeed, from the definition of $\mathcal P$ it follows that
${\mathcal P}^{\perp J_{1}}$ consists of exactly such vectors \[
z=\left[\begin{matrix}z_1\\ z_2\end{matrix}\right], \quad z_1\in
\left(B^k_{\Lambda^{\prime}\cap S}\right)^{m\times 1},\ z_2\in
\left(B^k_{\Lambda^{\prime}\cap S}\right)^{p\times 1}\] that \[
0=\scal{z_1,(T(A)^{}_{S,\Lambda^{\prime}})^{*}y}-\scal{z_2,(T(B)^{}_{S,\Lambda^{\prime}})^{*}y}=
\scal{T(A)^{}_{S,\Lambda^{\prime}}z_1-T(B)^{}_{S,\Lambda^{\prime}}z_2,y}\]
for all $y\in\left(B^k_{\Lambda^{\prime}\cap S}\right)^{p \times
1}$. In other words, $z\in{\mathcal P}^{\perp J_{1}}$ if and only if
\[
T(A)^{}_{S,\Lambda^{\prime}}z_1-T(B)^{}_{S,\Lambda^{\prime}}z_2=0.\]
This of course is equivalent to (\ref{pj1}).

By Lemma \ref{L:kerran} we have the alternative representation
$$
{\mathcal P}^{\perp J_{1}} = \text{Im } \begin{bmatrix}
T(D)^{}_{S,\Lambda^{\prime}} \\[2mm] T(C)^{}_{S,\Lambda^{\prime}}
\end{bmatrix}.
$$
In terms of this representation, the fact that ${\mathcal P}^{\perp
J_{1}}$ is a regular subspace means that the operator
$$
(T(D)^{}_{S,\Lambda^{\prime}})^{*}T(D)^{}_{S,\Lambda^{\prime}}
 - (T(C)^{}_{S,\Lambda^{\prime}})^{*} T(C)^{}_{S,\Lambda^{\prime}}$$
is invertible. By Theorem~\ref{th:inv}, the matrix function $(D^{*}D
- C^{*} C)^T$ admits a canonical factorization. Since all the values
of this matrix function are hermitian, its factorization can be
chosen in a special form described in Theorem \ref{nov206}. Passing
to transposed matrices, it can be written in the form
$$ D^{*}D - C^{*} C = R^{*} J_{0}R
$$
for an appropriate signature matrix $J_{0}$, where $R^{\pm 1} \in
(APW^k_{\Lambda^{\prime} \cap S})^{m \times m}$.  We then let
\begin{equation}\label{apr32}
\Theta = \begin{bmatrix} D \\ C \end{bmatrix} R^{-1}.
\end{equation}
One verifies  that $\Theta^{*} J_{1} \Theta=J_0. $ 
Thus $\Theta$ is $(J_{0}, J_{1})$-isometric.
Decompose $\Theta= \begin{bmatrix} \Theta_{11} & \Theta_{12} \\
\Theta_{21} & \Theta_{22} \end{bmatrix}, $ with the size of
$\Theta_{11}$ equal to $m \times (m-p)$ and that of $\Theta_{22}$
equal to $p \times p$. Repeating the arguments from the proofs of
\cite[Theorem 5.2]{BKRS} and \cite[Theorem 4.1]{BallRS}, we obtain
that for every $G \in (APW^k_{\Lambda^{\prime} \cap S})^{(m-p)
\times p}$, the linear fractional parametrization formula
(\ref{apr31}) yields a solution $F\in (APW^k)^{m\times
p}_{\Lambda^{\prime} \cap S}$ of the problem
\begin{equation}\label{mar293} A
F = B \quad {\rm and}\quad  \|F\|_\infty\leq 1.\end{equation}
To prove
that all solutions $F\in (APW^k_{\Lambda^{\prime} \cap S})^{m\times
p}$ are obtained this way, observe that
\[ G = (\Theta_{11}^*F -\Theta_{21}^*)(\Theta_{12}^* F - \Theta_{22}^*)^{-1} \]
due to $(J_0,J_1)$-isometricity of $\Theta$. This concludes the
proof of Theorem \ref{T:genToepcor}. \proofend

Note that the description (\ref{apr31}) of all $F \in (APW^k_{S \cap
\Lambda^{\prime}})^{m \times p}$ such that (\ref{solution}) holds,
applies provided the operator (\ref{Toepnew}) is positive definite
and invertible. If the operator (\ref{Toepnew}) is merely positive
semidefinite, then there exists an $F$ in the infinity norm matrix
Besicovitch space, having its Bohr-Fourier in $\Lambda'\cap S$
and satisfying (\ref{solution}). This statement is
completely analogous to the corresponding part of Theorem 5.2 of
\cite{BKRS}, and can be obtained in the same way. We refer the
interested reader to \cite{BKRS} for more details.
\bigskip

{\bf Acknowledgments.} We thank J. T. Anderson, A. Brudnyi, and A.
A. Pankov for useful consultations, and J. A. Ball for discussions
concerning preliminary versions of the draft for the paper.

\providecommand{\bysame}{\leavevmode\hbox
to3em{\hrulefill}\thinspace}
\providecommand{\MR}{\relax\ifhmode\unskip\space\fi MR }
\providecommand{\MRhref}[2]{%
  \href{http://www.ams.org/mathscinet-getitem?mr=#1}{#2}
} \providecommand{\href}[2]{#2}

\end{document}